\documentstyle[12pt]{article}
\newtheorem{defn}{Definition}
\begin{document}
\title{ Translatable radii of an operator in the direction of another operator II}
\author{Kallol Paul}
\date{}
\maketitle
\begin{abstract}
One of the couple of translatable radii of an operator in the direction of another operator
introduced in earlier work[13] is studied in details. A necessary and sufficient condition for a unit 
vector f to be a stationary vector of the generalized eigenvalue problem 
$ Tf = \lambda A f $ is obtained. Finally a theorem of Williams[16] is
generalized to obtain a translatable radius of an operator in the direction of another operator. 
\end{abstract}
\section{Introduction.}
Let T and A be two bounded linear operators on a complex Hilbert space H with
inner product (,) and norm $ \| \| $. 
Consider the generalized eigenvalue problem Tf = $ \lambda $ Af 
where f $ \in $ H and $ \lambda \in $ C, $ \lambda $ is called 
the eigenvalue of the above equation and f the corresponding eigenvector.
The non-negative functional 
\[ M_T(f) =  \|Tf - \frac{(Tf,Af)}{(Af,Af)} Af \|,~provided~\|Af\| \neq 0,  \]
gives the deviation of a unit vector f from being an eigenvector and 
\[ M_T(A) = sup_{\|f\|=1}~ \{ \|Tf - \frac{(Tf,Af)}{(Af,Af)} Af \| \},~provided~0 \notin  \sigma_{app}A,  \]
gives the supremum of all those deviations, where $ \sigma_{app}A $ is the set of approximate eigenvalues of A.\\ Geometrically $ Tf - \frac{(Tf,Af)}{(Af,Af)} Af $ is
 the component of Tf perpendicular to Af.
 For A = I problems related to the concepts considered here have been studied by  Bjorck and Thomee[2], Garske[8], 
 Prasanna[14], Fujii and Prasanna[6], Furuta et al[7], Fujii and Nakamoto[5], Izumino[9],  Nakamoto and Sheth[11],
 Mustafaev and Shul'man[10]  and many others. 
 \vspace{.1cm}
 {\bf \hrule}
\vspace{.2cm}
\noindent{\bf Keywords :} Stationary distance vectors, Translatable radii.
\pagebreak

\noindent Bjorck and Thomee[2] have shown that for a normal operator T,
 \[ M_T = sup_{\|f\|=1}~ \{ \|Tf - {(Tf,f)} f \| = R_T,  \] where $ R_T $ is the radius of 
 the smallest circle containing the spectrum. Garske[8] improved on the result to prove that for
  any bounded linear operator T,
  \[ M_T = sup_{\|f\|=1}~ \{ \|Tf - {(Tf,f)} f \| \geq R_T.  \]
 Stampfli[15] proved that for a bounded linear operator T $ \exists $ a unique complex scalar $ c_T $, defined as the center of mass of T such that 
 \[ \| T - c_{T}I \|^{2} + | \lambda |^{2} \leq \| T - c_{T}I + \lambda I \|^{2},~~ \forall ~ \lambda~ \in C. \] 
 
With the help of Stampfli's result Prasanna[14] proved that $ M_T = \| T - c_{T} I \| $. Later Fujii and Prasanna[6] improved on the inequality of Garske to show that $ M_{T} \geq w_{T} $ where $ w_{T} $ is the radius of the smallest circle containing the numerical range.\\  
\noindent In [12] we proved that for any two bounded linear operators 
T and A if $ 0 \notin \sigma_{app} A $ then there exists a unique complex scalar $ \lambda _0 $ such that 
$ \| T - \lambda_0 A \| \leq \|T - \lambda A \| ~\forall \lambda ~\in $ C.
We defined $ T - \lambda_0 A $ as the {\bf minimal-norm translation of T in 
the direction of A} and proved that $ \| T - \lambda_0 A \| = M_T(A) $. The equality of 
$ \inf_{\lambda}\| T - \lambda A \| = M_T(A) $ was also studied by E.Asplund and V.Pt$\acute{a}$k[1]\\
Then in [13] we introduced a couple of {\bf translatable radii of an operator T
in the direction of another operator A} as follows :\\
If 0 does not belong to the approximate point spectrum of  A let
\[ M_T(A) = sup_{\|f\|=1}~ \{ \|Tf - \frac{(Tf,Af)}{(Af,Af)} Af \| \} \]
\[ i.e., M_T(A) = sup_{\|f\|=1}~ {\{ {\| Tf\|}^{2} - \frac { {\mid(Tf,Af)\mid}^{2}} {(Af,Af)} \}}^{1/2} 
~~~~\]
and if $ 0  \notin \overline { W(A) }$, where $\overline { W(A) }$ stands for the closure of the numerical range of A, let 
\[ \tilde {M}_T(A) = sup_{\|f\|=1}~\{ \|Tf - \frac{(Tf,f)}{(Af,f)} Af \| \}~.~\]
We defined $ M_T(A) $ and $ \tilde {M}_T(A) $ as translatable radii of the operator
T in the direction of A and proved in [13] that 
if $ 0 \notin \overline {W(A)} $ then
\[ \tilde {M}_T(A) \geq M_{T}(A) \geq m_{T}(A)/\|A^{-1}\|~, \]
where $ m_{T}(A) $ is the radius of the smallest circle containing the set
$ W_{T}(A)= \{ ~(Tf,Af)/(Af,Af) ~: \|f\|=1 ~\} ~.$\\

\noindent Das[4] introduced the concept of stationary distance vectors 
while studying the eigenvalue problem Tf = $ \lambda $ f. 
Following the ideas of Das we here use the concept of 
stationary distance vectors to study the generalized eigenvalue problem 
Tf = $ \lambda $ Af and the  translatable radius $ M_{T}(A) $. We 
investigate the structure of the vectors for which the translatable 
radius $ M_T(A) $ is attained and prove that if 
$ M_{T}(A) $is attained at a vector f then $ M_{T^{*}}(A^{*}) $ is attained 
at the vector $ h/\|h\|$, where $ h =  Tf - (Tf,Af)/(Af,Af)~ Af  $ .
We also show that if g is a state ( normalized positive functional ) on 
 the Banach algebra B(H,H) of all bounded linear operators on H then 
\[ M_{T}(A) = sup \{~~ g(T^{*}T) - \frac{{\mid g(A^{*}T) \mid}^{2}}{g(A^{*}A)}~~ :
g ~~ is ~ a ~ state ~and~~g(A^{*}A) \neq 0 \} ~.\]
The last result mentioned here is a generalization of a theorem of Williams [16].
\section{Stationary distance vectors of the  generalized eigenvalue problem
Tf = $ \lambda $ Af}
In this section we study the following :\\`` For any two bounded linear operators T and A what are the vectors that are nearest to or
farthest from being eigenvectors of the equation $ Tf = \lambda Af $ in the sense that $ \|Tf - (Tf,Af)/(Af,Af)~Af \| $ 
with unit f is minimum or maximum ?''\\
We give a necessary and sufficient condition that a unit vector f is at a stationary distance 
from being an eigenvector. We call such f's the stationary distance vectors and the corresponding 
$ \lambda = $ (Tf,Af)/(Af,Af) the stationary distance value of the eigenvalue problem $ Tf = \lambda Af $.
We use the concept of stationary vectors the definition of which is given below:\\
\begin{defn}
{\bf Stationary vector.}\\
Let $ \varphi $ be a functional defined on the unit sphere of $ H $. Then a unit vector $f$  is said to be a stationary vector and $ \varphi $ is said to have a stationary value  at $f$ of $ \varphi $ iff the
function $ w_{g}(t) $ of a real variable t, defined as
\[ w_{g}(t) = \varphi ( \frac{f+tg}{\|f+tg\|}) \]
has a stationary value at t=0  i.e., $ w'_g(0) = 0 $ for any arbitrary but fixed vector $g$
$ \in $ H. \\
e.g., If $ \varphi(f) = \| Tf - (Tf,Af)/(Af,Af)~ Af \|^2 $ then  a stationary vector $f$ of functional $ \varphi $ is called the stationary distance vector of the eigenvalue problem $ Tf = \lambda Af $.
\end{defn}
We assume that  0 does not belong to the approximate point spectrum of A
 and prove the following theorem :

\noindent{\bf Theorem 1.} The necessary and sufficient condition for a unit vector f to be  a 
stationary distance vector of the generalized eigenvalue problem $ Tf = \lambda Af $ is that it satisfies the 
following 
\[ (T^{*}- \bar{\lambda}A^{*})(T- \lambda A) f = { \|h\|}^{2}f \]
where $ h =  Tf - \lambda Af $ and $ \lambda  = \frac{(Tf,Af)}{  (Af,Af)} $.\\

\noindent{\bf Proof.} Consider $ M_{T}(f) = \| Tf - (Tf,Af)/(Af,Af)~ Af \| $.
 Define the function $ w_{g}(t) $ of a real variable t as follows
\begin{eqnarray*}
w_{g}(t)~ =~ {M_T}^{2}~(\frac{f+tg}{\|f+tg\|})~=~ \frac{ {\|T(f+tg)\|}^{2} } { {\|f+tg\|}^{2} }~ -~
\frac{ {\mid ( T(f+tg), A(f+tg) ) \mid }^{~2} } {~  (A(f+tg), A(f+tg))~ {\|f+tg\|}^{2} } 
\end{eqnarray*}
where g is arbitrary but fixed vector in H.\\
At a stationary vector f we have $ w_{g}'(0) = 0 $ and so 
\begin{eqnarray*}
2~ Re~(T^{*}Tf,g) - {\|Tf\|}^{2}~ 2~ Re(f,g)~ -~ \frac{ {\|Af\|}^{2}}{{\|Af\|}^{4}}~
~[~ (Tf,Af)~~~~~~~~\\ 
 \{~ \overline { (Tf,Ag) + (Tg,Af) }~ \}~ + \overline{(Tf,Af)} 
 \{~(Tf,Ag)~ +~ (Tg,Af)~ \}~~ ] \\
 +~ \frac{ {\mid(Tf,Af)\mid}^{2}}{{\|Af\|}^{4}}~~ \{~
{\|Af\|}^{2}~ 2~ Re~(f,g)~ + ~2~ Re~ (A^{*}Af,g)~ \} = 0 ~.
\end{eqnarray*}
Since g is arbitrary we get,
\begin{eqnarray*}
 T^{*}Tf~ -~ {\|Tf\|}^{2}f~ - ~\lambda ~ T^{*}Af~ -~ \bar{\lambda} ~ A^{*}Tf~ +~ {\|Af\|}^{2}~
{\lambda}^{2}f + {\lambda}^{2} ~ A^{*}Af = 0 ~,\\
~~ where~~~ \lambda~ =~ (Tf,Af)/(Af,Af)~.~~~~~~~~~~~~~~~~~~~~~~~~~~~~~~~~~~~~~~~~~~
\end{eqnarray*}
Let $ h =  Tf - \lambda Af $, then  $ (h,Af) = 0 $ and  $ {\|h\|}^{2} =
{\|Tf\|}^{2} -  {\mid (Tf,Af) \mid }^{2} / (Af,Af)  $ .\\
So we get \[ (T^{*}- \bar{\lambda}A^{*})(T- \lambda A) f = { \|h\|}^{2}f~. \]
Thus the theorem is proved.\\

\noindent We now prove the following corollary :\\

\noindent{\bf Corollary 1.} If $ M_{T}(A) $ is attained at f then $ M_{T^{*}}(A^{*}) $ is also attained at 
$ h/\|h\|~~where ~~ h = Tf - (Tf,Af)/(Af,Af)~Af $.\\
{\bf Proof.} Suppose $ M_{T}(A) $ is attained at a vector f  and $ \lambda  = \frac{(Tf,Af)}{  (Af,Af)} $. Then f is a stationary distance vector
and so we get
\begin{eqnarray*}
& & (T^{*} - \bar{\lambda} A^{*} ) ( T - \lambda A) f = { \|h\|}^{2} f \\
\Rightarrow & & (T^{*} - \bar{\lambda} A^{*} ) h = { \|h\|}^{2} f \\
\Rightarrow & & (T^{*}h,A^{*}h) = \bar{\lambda} (A^{*}h,A^{*}h)\\
\Rightarrow & & \bar{\lambda} = \frac{(T^{*}h,A^{*}h)}{(A^{*}h,A^{*}h)} 
\end{eqnarray*}
\begin{eqnarray*}
Now~~ T^{*}h & = & \bar{\lambda}A^{*}h + {\|h\|}^{2} f \\
\Rightarrow  {\|T^{*}h\|}^{2} & = & {\mid \bar{\lambda} \mid}^{2} {\|A^{*}h\|}^{2}  + {\|h\|}^{4}  \\
\Rightarrow  {\|T^{*}h\|}^{2} & = & {\|h\|}^{2}  \{ 
{ \|Tf \|}^{2} - \frac { { \mid (Tf,Af) \mid }^{2}}{ (Af,Af) }\} + 
\frac { { \mid (Tf,Af) \mid }^{2}}{ (Af,Af) } ~.~ \frac{ {\|A^{*}h\|}^{2} }{{\|Af\|}^{2}} 
\end{eqnarray*}
If the minimal-norm translation of T in the direction of A is T itself then the minimal-norm translation
of $ T^{*} $ in the direction of $ A^{*} $ is also $ T^{*} $. So if $ M_{T}(A) = \|T\| $ 
then $~~ M_{T^{*}}(A^{*}) = \| T^{*}\| $.\\
Let $ M_{T}(A) = \|T\| = \|Tf\|, ~~ (Tf,Af)/(Af,Af) = 0 $ .\\
Then $ M_{T^{*}}(A^{*}) = \| T^{*}\| = \|T\|= \|T^{*}h\|/\|h\| ,~~ since~~
(Tf,Af)/(Af,Af) = 0 $.\\
This completes the proof.\\

\noindent Next we prove the following theorem :\\

\noindent{\bf Theorem 2.} Suppose T and A are two selfadjoint operators and  f be a unit stationary
distance vector such that (Tf,Af)
is real, then f can be expressed as the linear combination of two eigenvectors 
of the problem $ Tf = \lambda Af $.\\
{\bf Proof.} As both T and A are selfadjoint and f is a stationary distance 
vector with (Tf,Af) real we get from the last theorem
\[ {(T- \lambda A)}^{2} f = {\|h\|}^{2} f .\] So we get
\begin{eqnarray*}
\Rightarrow  {(T- \lambda A)}^{2} f \pm \|h\|h &=& {\|h\|}^{2} f \pm \|h\|h \\
\Rightarrow T(Tf - \lambda Af \pm \|h\|f) &=& (\lambda A \pm \|h\|)(Tf - \lambda Af \pm \|h\|f)
\end{eqnarray*}
\[Let ~~ g_{1} = Tf - \lambda Af + \|h\|f \]
\[and ~~ g_{2} = Tf - \lambda Af - \|h\|f~.\]
Then we get
\[ Tg_{1} = (\lambda A + \|h\|)g_{1} ~~ and ~~ Tg_{2} = (\lambda A - \|h\|)g_{2} \]
so that 
\[ (T - \lambda A)g_{1} = \|h\|g_{1} ~~ and ~~ (T - \lambda A)g_{2} = - \|h\|g_{2} ~.\]
Thus $ f = (g_{1} - g_{2} )/ (2\|h\|) $ completes the proof.

\section{On the attainment of $ M_{T}(A)$ }
Suppose $ \{ f_{n} \} $ be a sequence of unit vectors such that
\[ { \|Tf_{n} \|}^{2} - \frac { { \mid (Tf_{n},Af_{n}) \mid }^{2}}{ (Af_{n},Af_{n}) } 
~~ \longrightarrow ~~{M_{T}(A)}^{2}  .\]
As the unit sphere in H is weakly compact without loss of generality we may assume that $ \{f_{n}\} $ converges weakly
to f i,e, $ f_{n} \rightharpoonup  f $.\\

\noindent We now prove the following theorem :\\

\noindent{\bf Theorem 3.} Suppose $ \{ f_{n}\} $ be a weakly convergent sequence of unit vectors such that 
\[ { \|Tf_{n} \|}^{2} - \frac { { \mid (Tf_{n},Af_{n}) \mid }^{2}}{ (Af_{n},Af_{n}) } 
~~ \longrightarrow ~~{M_{T}(A)}^{2}  .\]
If the weak limit f is non-zero then $ M_{T}(A) $ is attained for the vector $ 
f/\|f\| $. If the supremum is not attained then all such sequences
must tend weakly to zero.\\
{\bf Proof.} Since $ M_{T}(A) $ is translation invariant in the direction of A so without any loss of generality we may
assume that the minimal-norm translation of T in the direction of A is T itself
i,e, $ M_{T}(A) = \|T\| $ . \\
So there exists a sequence  $ \{ f_{n} \},~ f_{n} \in H , ~\|f_{n}\| =1 ~ $ such that 
$ \| Tf_{n} \| \longrightarrow \|T\| $ and $ (Tf_{n}, Af_{n}) \longrightarrow 0 ~$.
Considering the positive operator $ {\|T\|}^{2}I - T^{*}T $ we have 
\begin{eqnarray*}
& & ({\|T\|}^{2}f_{n} - T^{*}T f_{n},f_{n} ) \longrightarrow 0 \\
\Rightarrow & & {\|T\|}^{2}f_{n} - T^{*}T f_{n} \longrightarrow 0 ~, ~~by ~
property~ of ~ positive ~ operators.\\
If ~ f &\neq& 0 ~we~have~\\
& & {\|T\|}^{2}(f_{n},f) - (T^{*}T f_{n},f) \longrightarrow 0 ~.
\end{eqnarray*}
Since $ f_{n} \rightharpoonup f $ and weak limit f is unique we get
\[ {\|T\|}^{2} = \frac{{\|Tf\|}^{2}}{{\|f\|}^{2}} ~.\]
The result that `` if $ f_n \rightharpoonup f ,~  \|Tf_n\| \rightarrow \|T\|  $ 
and $ f \neq 0 $ then $ \|T\| $ is attained at $ f/\|f\| $ '' 
follows directly from the corollary 1 of Das[3].\\
As $ M_T(T) = \|A\| $ the theorem is proved.

\section{ On generalization of a Theorem of Williams }
Let $ {\cal B} $ denote the set of all normalized positive linear functionals (states) on B(H,H) i.e.,
\[ {\cal B} = \{ ~g ~: ~ g\in L(B(H,H),C) \hspace{2mm} and \hspace{2mm} g(I) = 1 = \|g\| ~\}  \]
Clearly $ {\cal B} $ is  $ weak^{*} $compact. 
Let $ {\cal P} = \{ ~ g~: ~g \in {\cal B} ~ and ~ g(A^{*}A) \neq 0 ~\} $ .\\
Williams[16] proved that for any bounded linear operator T, $~ \|T\| \leq 
\|T - \lambda I \| ~~ \forall \lambda ~\in $ C  iff there exists a 
state f such that $ f( T^*T) = \|T^*T\| $ and f(T) = 0. We here show that if 
for two bounded linear operators T and A, $~ \|T \| \leq \| T -\lambda A \| ~~
\forall ~\lambda ~\in $ C then $ \| T\|^2~ = sup \{~~ g(T^{*}T) - 
\frac{{\mid g(A^{*}T) \mid}^{2}}{g(A^{*}A)}~~ :g ~~ is ~ a ~ state ~ 
and ~~ g(A^{*}A) \neq 0 \} ~.$\\
We now prove the following theorem :\\

\noindent{\bf Theorem 4.} $ [M_{T}(A)]^2 = sup \{~~ g(T^{*}T) - \frac{{\mid g(A^{*}T) \mid}^{2}}{g(A^{*}A)}~~ :
g ~~ is ~ a ~ state ~ and ~~ g(A^{*}A) \neq 0 \} ~.$\\
{\bf Proof.} Let $ [S_{T}(A)]^2 = sup \{~~ g(T^{*}T) - \frac{{\mid g(A^{*}T) \mid}^{2}}{g(A^{*}A)}~~ :
g ~~ is ~~ a ~~ state ~~ and~~g(A^*A) \neq 0 \} $. \\
Clearly $ S_{T + \lambda A}( A) = S_{T}(A ) $ and  $ M_{T + \lambda A}( A) = M_{T}(A ) $ 
so that both  are translation invariant in the direction of A. 
 Without loss of generality we assume that $ M_{T}(A) = \|T\| ~.$\\
Now for each  $ x \in H $,  $ \|x\|= 1 $, let $ g_{x} : B(H,H) \longrightarrow C $
be defined as $ g_{x}(U) = (Ux,x)~~ \forall U~ \in~ B(H,H) $.\\
Then $ g_{x} $ is a state and $ g_{x}(A^{*}A) \neq 0 $.\\
So
\begin{eqnarray*}
\|T\|  & = & sup_{g_{x}} { \{ g_{x}(T^{*}T) - \frac{{\mid g_{x}(A^{*}T) \mid}^{2}}
{g_{x}(A^{*}A)} \} }^{1/2} \\
         & \leq & sup_{g \in {\cal P} } { \{ g(T^{*}T) - \frac{{\mid g(A^{*}T) \mid}^{2}}
  {g(A^{*}A)} \} }^{1/2}  \\
         & \leq & sup_{g \in {\cal P} }  { \{g(T^{*}T) \} }^{1/2} \\
	 & = & \|T\|~.
\end{eqnarray*}	 
This completes the proof.\\

\noindent {\bf Note.} For A=I the result of Williams follows easily from Theorem 4.\\

\noindent {\bf Acknowledgement.} The author thanks Professor T.K.Mukherjee and Professor K.C.das
 for their help while preparing this paper. The author would also like to thank the referee for
  his invaluable suggestion.\\

\vspace{.02cm}
{\bf \hrule }
\vspace {.1cm} 
Reader in Mathematics\\
Department of Mathematics\\
Jadavpur University\\
Kolkata 700032\\
INDIA.\\
$ e-mail~: ~kalloldada@yahoo.co.in , ~~~kpaul@math.jdvu.ac.in$ \\

\begin{thebibliography}{99}
\bibitem{As:Ac}E.Asplund and V.Pt$\acute{a}$k, A minimax inequality for operators and a related numerical range,
{\em Acta Mathematica}, 126 (1971), 53-62.
\bibitem{Bj:Ar}G.Bjorck and V.Thomee, A property of bounded normal
operators in Hilbert Space, {\em Arkiv for Math.}, 4 (1963), 551-555.
\bibitem{Da:An} K.C.Das, Extrema of the Rayleigh Quotient and normal Behavior
of an operator, {\em Journal of Mathematical Analysis and Applications}, Vol.41
No.3 (1973) 765-774.
\bibitem{Da:Sc} K.C.Das, Stationary distance vectors and their relation with eigenvectors,
{\em Science Academy Medals for Young Scientists-Lectures}, (1978) 44-52.
\bibitem{Fu:Ma} M.Fujii and R. Nakamoto, An estimation of the transcendental radius of an operator,
{\em Math. Japonica}, 27 (1982), 637-638.
\bibitem{Fu:Ja} M.Fujii and S.Prasanna, Translatable radii for operators, {\em Mathematica
Japonica}, 26 (1981) 653-657.
\bibitem{Fu:Ma} T.Furuta, S.Izumino and S.Prasanna, A characterisation of centroid operators,
{\em Math. Japonica}, 27 (1982) 105-106.
\bibitem{Ga:Pr}G.Garske, An equality concerning the smallest disc that 
contains the spectrum of an operator, {\em Proc. Amer. Math. Soc.}, 78
(1980), 529-532.
\bibitem{Iz:Ma} S.Izumino, An estimation of the transcendental radius of an operator,
{\em Math. Japonica}, 27 No.5 (1982), 645-646.
\bibitem{Mu:Ma} G.S.Mustafaev and V.S.Shul'man, An estimate of the norms of inner derivation  in some operator algebras. {\em Math. Notes(English. Russian original)}
 45, No.4 (1989) 337-341;  translation from Mat. Zametki 45, No.4, 105-110 (1989).
\bibitem{Na:Ma} R. Nakamoto and I.H.Sheth, On centroid operators. {\em Math. Japonica},
 29, No.2 (1984) 287-289. 
\bibitem{Pa:Pu} K.Paul, Sk.M.Hossein and K.C.Das, Orthogonality on B(H,H) and minimal-norm 
operator, {\em Journal of Analysis and Applications}, Vol. 6, No. 3 (2008) 169-178.
\bibitem{Pa:Sc} K.Paul, Translatable radii of an operator in the direction of
another operator, {\em Scientae Mathematicae}, Vol.2 No.1 (1999) 119-122.
\bibitem{Pr:Ja}S.Prasanna, The norm of a derivation and the Bjorck-Thomee-Istratescu 
theorem, {\em Mathematica Japonica}, 26 (1981), 585-588.
\bibitem{St:Pa} G. Stampfli, The norm of a derivation, {\em Pacific J. math.}, 33 (1970) 737-747.
\bibitem{Wi:Pr} J.P.Williams, Finite operators, {\em Proc. Amer. Math. Soc.}
Vol.26 (1970) 129-136.
\end{thebibliography}
\end{document}